\documentclass{siamltex}
\usepackage{amssymb,amsmath,amsfonts}
\usepackage{graphics}
\def\exp{\mathop{\mathrm{exp}}}
\def\trace{\mathop{\mathrm{trace}}}
\def\log{\mathop{\mathrm{log}}}
\def\md{M_D}
\def\moff{M_{\mathrm{off}}}

\def\lmin{\lambda_{\mathrm{min}}}
\newtheorem{example}{Example}

\begin{document}
\title{Determinant Approximations}

\author{Ilse C.F. Ipsen\thanks{%
Center for Research in Scientific Computation,
Department of Mathematics, North Carolina State University, P.O. Box 8205,
Raleigh, NC 27695-8205, USA ({\tt ipsen@math.ncsu.edu},
{\tt http://www4.ncsu.edu/{\char'176}ipsen/}).
Research supported in part by  NSF grants DMS-0209931 and DMS-0209695.}
\and
Dean J. Lee\thanks{%
Department of Physics, North Carolina State University, Box 8202,
Raleigh, NC 27695-8202, USA ({\tt djlee3@unity.ncsu.edu},
{\tt http://www4.ncsu.edu/{\char'176}djlee3/})
Research supported in part by NSF grant DMS-0209931.}
}
\maketitle

\begin{abstract}
A sequence of approximations for the determinant 
and its logarithm of a complex matrix
is derived, along with relative error bounds.  The determinant
approximations are derived from expansions of
$\det(X)=\exp(\trace(\log(X)))$, and they apply to
non-Hermitian matrices.  Examples illustrate that these determinant  
approximations are efficient
for lattice simulations of finite temperature nuclear matter, and
that they use significantly less space than Gaussian elimination.

The first approximation in the sequence is
a block diagonal approximation; it represents an extension of
Fischer's and Hadamard's inequalities to non-Hermitian
matrices.  In the special case of Hermitian positive-definite
matrices, block diagonal approximations can be
competitive with sparse inverse approximations.
At last, a different representation of sparse inverse approximations
is given and it is shown that their accuracy increases as more
matrix elements are included.
\end{abstract}



\section{Introduction}\label{s_intro}
For a complex matrix we present approximations for the determinant
and its logarithm, together with error bounds.  

The approximations were motivated by a problem in computational
quantum field theory: the simulation of finite temperature nuclear
matter on a lattice \cite{LeeI03}.  In this application, the logarithm
of the determinant is desired to  2-3 significant digits. The
matrices are sparse, and non-Hermitian.
Because the desired accuracy is low, LU decomposition with partial
pivoting \cite[\S 14.6]{Hig02}, \cite[\S 3.18]{Wil63} is too costly.
Since the matrices are not Hermitian positive-definite, sparse
approximate inverses \cite{Reu02}, Gaussian quadrature based 
methods \cite{BaiG97}, and Monte Carlo methods \cite[\S 4]{Reu02} 
or hybrid Monte Carlo methods \cite{DuKPR87,GoLTRS87,ScSS86} 
do not apply. Monte Carlo and quadrature-based methods 
can be extended to non-Hermitian matrices, however then 
the sign of the determinant is usually lost, e.g. \cite[\S 3.2.3]{BaiFG96}.

To approximate the determinant $\det(M)$ we decompose $M=\md+\moff$
such that $\md$ is a non-singular matrix. Then
$\det(M)=\det(\md)\det(I+\md^{-1}\moff)$, where $I$ is the identity matrix.
In
$$\det(I+\md^{-1}\moff)=\exp(\trace(\log(I+\md^{-1}\moff))),$$
we expand $\log(I+\md^{-1}\moff)$, obtaining
a sequence of increasingly accurate approximations.
Error bounds for these approximations depend on the spectral
radius of $\md^{-1}\moff$.

\subsection*{Overview}
The determinant approximations and their
error bounds are presented in \S \ref{s_main}.
Approximations from block diagonals (\S \ref{s_diag}) are extended to 
a sequence of  higher order approximations (\S \ref{s_higher}).
They simplify for checkerboard matrices (\S \ref{s_tridiag})
which occur in the neutron matter simulations in \cite{LeeI03}.
Comparisons with sparse inverse
approximations of determinants, which are limited to Hermitian 
positive-definite matrices (\S \ref{s_si})
illustrate the competitiveness of block diagonal approximations.
As expected, the
accuracy of sparse inverse approximations increases as more matrix
elements are included.
Numerical results with matrices from nuclear matter simulations
(\S \ref{s_app}) show that determinant approximations of desired accuracy
can be obtained fast, in 1-3 iterations; and that they require
significantly less space than Gaussian elimination (with partial
or complete pivoting).

\subsection*{Notation}
The eigenvalues of a complex square matrix $A$ are $\lambda_j(A)$, and the
spectral radius is $\rho(A)\equiv\max_j{|\lambda_j(A)|}$.
The identity matrix is $I$, and $A^*$ is the conjugate transpose of $A$.
We denote by $\log(X)$ and $\exp(X)$ the logarithm and exponential 
function of a matrix $X$, and by $\ln(x)$ and $e^x$
the natural logarithm and exponential function of a scalar $x$.

\section{Determinant Approximations}\label{s_main}
We present approximations to the determinant and its logarithm, 
as well as error bounds.  

\subsection{Diagonal Approximations}\label{s_diag}
We present relative error bounds for the approximation 
of the determinant by the determinant of a block diagonal.
Let $M$ be a complex square matrix of order $n$
partitioned as a $k\times k$ block matrix
$$M=\begin{pmatrix}M_{11}   & M_{12} &\ldots & M_{1k} \cr
             M_{21}   & M_{22} &\ldots & M_{2k} \cr
             \vdots   & \ddots &\ddots & \vdots\cr
             M_{k1}   & M_{k2} &\ldots & M_{kk}\end{pmatrix},$$
where the diagonal blocks $M_{jj}$ are square but not necessarily of the 
same dimension. Analogously,
decompose $M=\md+\moff$ into diagonal blocks $\md$ and off-diagonal
blocks $\moff$,
\begin{eqnarray}\label{e_part}
\md=\begin{pmatrix}M_{11} &  &       &    \cr
                                  & M_{22} &  &     \cr
                                 &    &\ddots & \cr
                                 &     &       &  & M_{kk}\end{pmatrix},\qquad
\moff=\begin{pmatrix}0        & M_{12} &\ldots & M_{1k} \cr
             M_{21}   & 0      &\ldots & M_{2k} \cr
             \vdots   & \ddots &\ddots & \vdots\cr
             M_{k1}   & M_{k2} &\ldots & 0\end{pmatrix}.
\end{eqnarray}
The block diagonal matrix $\md$ is called a pinching of $M$ 
\cite[\S II.5]{Bha97}.
We consider the approximation of $\det(M)$ by the determinant of a
pinching, $\det(\md)$; and in particular bounds 
of the form $\det(M)\leq \det(\md)$. The matrices for which such bounds
are known to hold are characterized by eigenvalue monotonicity
of the following kind.

A complex square matrix $M$ is a $\tau$-matrix if \cite[pp 156-57]{EngS76}:
\begin{enumerate}
\item Each principal submatrix of $M$ has at least one real eigenvalue.
\item If $S_1$ is a principal submatrix of $M$ and $S_{11}$ a principal
submatrix of $S_1$ then $\lmin(S_1)\leq \lmin(S_{11})$, where
$\lmin$  denotes the smallest \textit{real} eigenvalue.
\item $\lmin(M)\geq 0$.
\end{enumerate}
The class of $\tau$-matrices includes Hermitian positive-definite,
M-matrices and totally non-negative matrices \cite[pp 156-57]{EngS76},
\cite[Theorem 1]{Mehr84a}.

\paragraph{Hadamard-Fischer Inequality}
If $M$ is a $\tau$-matrix then \cite[Theorem 4.3]{EngS76}
\begin{eqnarray}\label{e_fh}
\det(M)\leq \det(\md).
\end{eqnarray}
Strictly speaking, (\ref{e_fh}) is called a Hadamard-Fischer inequality 
only for $k=2$ \cite[(0.5)]{EngS76}.
If $k=2$ and $M$ is Hermitian positive-definite then (\ref{e_fh}) 
is Fischer's inequality \cite[Theorem 7.8.3]{HoJ85}.
If $k=n$ and $M$ Hermitian positive-definite then (\ref{e_fh}) is Hadamard's 
inequality \cite[Theorem 8]{CovT88}, \cite[Theorem 7.8.1]{HoJ85}.
Extensions of (\ref{e_fh}) to generalized 
Fan inequalities are derived in \cite{Mehr84b,Mehr84a}.
The Hadamard-Fischer inequality (\ref{e_fh}) implies the obvious relative 
error bound for the determinant of a pinching,
$$0< {\det(\md)-\det(M)\over \det(\md)}\leq 1.$$
In the theorem below we tighten the upper bound.

\begin{theorem}\label{t_1a}
Let $M$ be a complex matrix of order $n$.
If $\det(M)$ is real, $\md$ is non-singular with $\det(\md)$ real,
and all eigenvalues $\lambda_j(\md^{-1}\moff)$ are real 
with $\lambda_j(\md^{-1}\moff)>-1$, then
$$0< {\det(\md)-\det(M)\over \det(\md)}\leq 
1-e^{-{n\rho^2\over 1+\lmin}},$$
where $\rho\equiv\rho(\md^{-1}\moff)$ and
$\lmin\equiv\min_{1\leq j\leq n}{\lambda_j(\md^{-1}\moff)}$.
\end{theorem}

\begin{proof}
Write $\det(M)=\det(\md)\det(I+A)$, where $A\equiv \md^{-1}\moff$.
Since $I+A$ is non-singular,
\cite[Theorem 6.4.15(a)]{HoJ91} and \cite[Problem 6.2.4]{HoJ91}
imply $\det(I+A)=\exp(\trace(\log(I+A)))$. Furthermore,
$\log(I+A)=\sum_{p=1}^{\infty}{{(-1)^{p-1}\over p}A^p}$
\cite[(7) in \S 9.8]{LaT85}.
Hence 
$$\det(I+A)=\exp\left(\sum_{j=1}^n{\ln(1+\lambda_j(A))}\right).$$
Because $\lambda_j(A)>-1$, $1\leq j\leq n$,  
we can apply the inequality 
${\lambda\over 1+\lambda}\leq\ln(1+\lambda)\leq\lambda$ \cite[4.1.33]{AbS72}
to obtain
$$\exp(\trace(A))\> e^{-{n\rho(A)^2\over 1+\lmin}}
\leq\det(I+A)\leq \exp(\trace(A)).$$
At last use the fact that  $\md$ is block diagonal and  
$\trace(\moff)=\trace(A)=0$.
\end{proof}

The upper bound for the relative error is small if the eigenvalues of
$\md^{-1}\moff$ are small in magnitude but not too close to $-1$. 
The pinching $\det(M_D)$ can be a bad approximation to $\det(M)$
when $I+\md^{-1}\moff$ is close to singular. 
If  $\det(\md)>0$ then Theorem \ref{t_1a} implies a lower bound 
for $\det(M)$,
$$e^{-{n\rho^2\over 1+\lmin}}\det(\md)\leq \det(M)\leq \det(\md).$$
In the argument of the exponential function in Theorem \ref{t_1a}
we have $\lmin<0$ because $\md^{-1}\moff$ has a zero diagonal so that 
$\trace(\md^{-1}\moff)=0$.  Hence $n\rho^2/(1+\lmin)>n\rho^2$.  

\begin{corollary}\label{c_1}
Theorem \ref{t_1a} holds for Hermitian positive-definite matrices.
In  particular, Theorem  \ref{t_1a} implies
error bounds for Fischer's and Hadamard's inequalities.
\end{corollary}

The following example shows that $|\det(M)|\leq|\det(\md)|$ may not hold
when $\md^{-1}\moff$ has complex eigenvalues, or real eigenvalues
smaller than $-1$.

\begin{example}
Even if all eigenvalues of $\md^{-1}\moff$ satisfy 
$|\lambda_j(\md^{-1}\moff)|<1$, it is still possible that
$|\det(M)|>|\det(\md)|$ if some eigenvalues are complex.

Consider
$$M=\begin{pmatrix}1&\alpha \cr \alpha &1\end{pmatrix},\qquad
\md=\begin{pmatrix}1&0\cr 0&1\end{pmatrix},\qquad
\moff=\begin{pmatrix}0&\alpha\cr \alpha &0\end{pmatrix}=\md^{-1}\moff.$$
Then $\lambda_j(\md^{-1}\moff)=\pm \alpha$, and $\det(M)=1-\alpha^2$.
Choose $\alpha ={1\over 2}\imath$, where $\imath=\sqrt{-1}$.
Then both eigenvalues of $\md^{-1}\moff$ are complex, and
$|\lambda_j(\md^{-1}\moff)|<1$. But $\det(M)=1.25>1=\det(\md)$.

The situation
$\det(\md)>\det(M)$ can also occur when $\md^{-1}\moff$ has a real eigenvalue
that is less than $-1$. If $\alpha=3$ in the matrices above then
one  eigenvalue of $\md^{-1}\moff$ is $- 2$, and $|\det(M)|=8>\det(\md)=1$. 
In general, 
$|\det(M)|/\det(\md)\rightarrow\infty$ as $|\alpha|\rightarrow\infty$.
\end{example}

This example illustrates that, unless the eigenvalues of $\md^{-1}\moff$
are real and greater than $-1$, $\det(\md)$ is, in general, not 
a bound for $\det(M)$. In the case of complex eigenvalues,
however, we can still determine how well $\det(\md)$ \textit{approximates}
$\det(M)$. Below is a relative error bound for the case when 
$M$ is 'diagonally dominant', in the sense
that the eigenvalues of $\md^{-1}\moff$ are small in magnitude.

\begin{theorem}[Complex Eigenvalues]\label{t_2}
Let $M$ be a complex matrix of order $n$.
If $M_D$ is non-singular and $\rho\equiv\rho(\md^{-1}\moff)<1$
then
$${|\det(M)-\det(\md)|\over |\det(\md)|}\leq c\rho\>e^{c\rho},\qquad
\mathrm{where}\quad c\equiv -n\ln(1-\rho).$$
If also $c\rho<1$ then
$${|\det(M)-\det(\md)|\over |\det(\md)|}\leq {7\over 4}c\rho.$$
\end{theorem}

\begin{proof} This is a special case of Theorem \ref{t_3}.
\end{proof}

\begin{corollary}\label{c_2}
Theorem \ref{t_2} holds for the following classes of matrices:
M-matrices; Hermitian positive-definite matrices if $k=n$; 
Hermitian positive definite block tridiagonal matrices with 
equally-sized blocks of dimension $n/k$.
\end{corollary}

\begin{proof}
In all cases $\rho(\md^{-1}\moff)<1$.
\end{proof}

In the special case of strictly diagonally dominant matrices, Theorem \ref{t_2}
leads to a bound for the approximation of $\det(M)$ by the product of
the diagonal elements.

\begin{corollary}\label{c_2a}
If the complex square matrix $M=(m_{ij})_{1\leq i,j\leq n}$ is 
strictly row diagonal dominant then
$${|\det(M)-\prod_{i=1}^n{m_{ii}}|\over |\prod_{i=1}^n{m_{ii}}|}
\leq c\rho\>e^{c\rho},
\qquad\mathrm{where}\quad
\rho\leq \max_i{\sum_{j=1,j\neq i}^n{\left|\frac{m_{ij}}{m_{ii}}\right|}},
\quad c\equiv -n\ln(1-\rho).$$
If also $c\rho<1$ then
$${|\det(M)-\prod_{i=1}^n{m_{ii}}|\over |\prod_{i=1}^ni{m_{ii}}|}
\leq {7\over 4}c\rho.$$
\end{corollary}

\begin{proof} This is a consequence of Gerschgorin's theorem
\cite[Theorem 7.2.1]{GovL96}.
\end{proof}

Corollary \ref{c_2a} implies that the product of diagonal elements 
is a good approximation for  $\det(M)$ if $M$ is strongly diagonally dominant.

\subsection{A Sequence of General Higher Order Approximations}\label{s_higher}
We extend the diagonal approximations in \S \ref{s_diag}
to a sequence of more general approximations
that become increasingly more accurate. These approximations, called 
'zone determinant approximations' in \cite{LeeI03}, are justified in the
context of nuclear matter simulations. 
As before, decompose $M=\md+\moff$ into diagonal blocks $\md$ and off-diagonal
blocks $\moff$ (actually, our results hold for any decomposition
$M=M_0+M_E$ where $M_0$ is non-singular and $\rho(M_0^{-1}M_E)<1$).
Below we give a sequence of approximations $\delta_m$ for $\ln(\det(M))$
and $\Delta_m$ for $\det(M)$, as well as absolute bounds for $\delta_m$
and relative bounds for $\Delta_m$.
An absolute bound for the logarithm suffices 
because $\ln(\det(M))>1$ in our applications.

\begin{theorem}\label{t_3}
Let $M=\md+\moff$ be a complex matrix of order $n$,
$\md$ be non-singular and $\rho\equiv\rho(\md^{-1}\moff)<1$. Define
$$\delta_m\equiv\ln(\det(\md))+\>\sum_{p=1}^m
{{(-1)^{p-1}\over p}\trace((\md^{-1}\moff)^p)},\qquad
\Delta_m\equiv e^{\delta_m},\qquad m\geq 1.$$
Then
$$|\ln(\det(M))-\delta_m|\leq c\rho^m,\qquad
{|\det(M)-\Delta_m|\over |\Delta_m|}\leq c\rho^m\>e^{c\rho^m}$$
where $c\equiv -n\ln(1-\rho)$.
If also $c\rho^m<1$ then
$${|\det(M)-\Delta_m|\over |\Delta_m|}\leq {7\over 4}\>c\>\rho^m.$$
\end{theorem}

\begin{proof}
As in the proof of Theorem \ref{t_1a} 
$\det(M)=\det(\md)\det(I+A)$, where $A\equiv \md^{-1}\moff$
and $\log(I+A)=\sum_{p=1}^{\infty}{{(-1)^{p-1}\over p}A^p}$.
Hence
$$\trace\left(\log(I+A)\right)=\sum_{p=1}^{\infty}
{{(-1)^{p-1}\over p}\trace(A^p)}.$$
Define the truncated sums
$$L_m\equiv\sum_{p=1}^m{{(-1)^{p-1}\over p}\trace(A^p)},
\qquad D_m\equiv e^{L_m},\qquad m\geq 1.$$
Then 
$$\trace\left(\log(I+A)\right)=L_m+z,\qquad
z\equiv \sum_{i=1}^n{\left\{\ln(1+\lambda_i(A))-
\sum_{p=1}^m{{(-1)^{p-1}\over p}\lambda_i(A)^p}\right\}}.$$
Applying to each of the $n$ terms the inequality
$$\left|\ln(1+\lambda)-\sum_{p=1}^m{{(-1)^{p-1}\over p}\lambda^p}\right|
\leq -\ln(1-|\lambda|)\>|\lambda|^m$$
\cite[4.1.24]{AbS72}, \cite[4.1.38]{AbS72} 
gives $|\ln(\det(I+A))-L_m|\leq c\rho^m$. The first bound follows
now with $\delta_m=\ln(\det(\md))+L_m$.

From the first bound, the fact that $\det(I+A)=D_me^z$, and
$|e^z-1|\leq |z|\> e^{|z|}$ \cite[4.2.39]{AbS72} follows
$${|\det(I+A)-D_m|\over |D_m|}\leq c\rho^m\>e^{c\rho^m}.$$
We get the second bound from $\Delta_m=\det(\md)\>D_m$.

If also $c\rho^m<1$ then \cite[4.2.38]{AbS72}
$${|\det(I+A)-D_m|\over |D_m|}\leq {7\over 4}c\>\rho^m.$$
\end{proof}

The accuracy of the approximations in Theorem \ref{t_3}
is determined by the spectral radius
$\rho$ of $\md^{-1}\moff$. In particular, the absolute error bound for the 
approximation  $\delta_m$ is proportional to $\rho^m$, hence
the approximations  tend to improve with increasing $m$.
The numerical results in Sections \ref{s_si} and \ref{s_app} illustrate 
that the pessimistic factor in the bound 
$|\ln(\det(M))-\delta_m|\leq -n\ln(1-\rho)\>\rho^m$ is $n$. 
We found that replacing $n$ by the number of eigenvalues whose magnitude 
is close to $\rho$ makes the bound tight.
The approximations for the logarithm can be determined from successive updates
$$\delta_0\equiv \ln(\det(\md)),\qquad\delta_m=\delta_{m-1}+
{(-1)^{m-1}\over m}\trace((\md^{-1}\moff)^m), \qquad m\geq 1,$$
and $\Delta_m=e^{\delta_m}$. Note that $e^{\delta_0}=\det(\md)$
is the block diagonal approximation from (\ref{e_part}). Hence 
Theorem \ref{t_2} is a special case of Theorem \ref{t_3} with $m=1$. 
If a block diagonal determinant approximation is sufficiently 
accurate, as in \S \ref{s_app}, it can be much cheaper 
to compute than a determinant via Gaussian elimination.

\subsection{Checkerboard Matrices}\label{s_tridiag}
For this particular class of matrices, which occurs in our 
applications \cite{LeeI03},
every other determinant approximation $\Delta_m$ has increased accuracy.
We call a matrix $M$ with equally sized blocks $M_{ij}$ of
dimension $n/k$ in (\ref{e_part}) an \textit{odd checkerboard matrix}
(with regard to the block size $n/k$)
if $M_{ij}=0$ for $i$ and $j$ both even or both odd, $1\leq i,j \leq k$;
and an \textit{even checkerboard matrix} if $M_{ij}=0$ for $i$ odd
and $j$ even or vice versa. An odd checkerboard matrix has zero 
diagonal blocks, hence its trace is zero.

\begin{theorem}\label{t_checker}
If, in addition to the conditions of Theorem \ref{t_3},
$\moff$ is an odd checkerboard matrix then
$$\delta_0=\ln(\det(\md)),\qquad \delta_m=
\begin{cases}\delta_{m-1} & \mathrm{if}~m~\mathrm{is~odd} \cr 
\delta_{m-2}-\trace\left((\md^{-1}\moff)^m\over m\right)
& \mathrm{if}~m~\mathrm{is~even}.\end{cases}$$
\end{theorem}

\begin{proof}
If $A$ and $B$ are odd checkerboard matrices (with regard to the
same block size) then $AB$  is an even checkerboard matrix. 
If $A$ is an odd checkerboard matrix and $B$ an even checkerboard
matrix then $AB$ and $BA$ are odd checkerboard matrices.

Since $\md^{-1}\moff$ is an odd checkerboard matrix, so are
the powers $(\md^{-1}\moff)^p$ for odd $p$. This means
$\trace\left((\md^{-1}\moff)^p\right)=0$ for odd $p$. 
Hence the approximations in Theorem \ref{t_3}
satisfy $\delta_{m}=\delta_{m-1}$ for $m$ odd. For $m$ even
$$\delta_m=\delta_{m-2}+\left({(-1)^{m-1}\over m}\>
\trace\left((\md^{-1}\moff)^m\right)\right)=
\delta_{m-2}-\left({\trace\left((\md^{-1}\moff)^m\right)\over m}\right).$$
\end{proof}

Theorem \ref{t_checker}  shows that an odd-order approximation is
equal to the previous even-order approximation. Hence the even-order
approximations gain one order of accuracy.

\section{Comparison with Sparse Inverse Approximations}\label{s_si}
In the special case of Hermitian positive-definite matrices,
we illustrate that block-diagonal determinant approximations 
(see Corollary \ref{c_1}) can compare favourably with approximations
based on sparse approximate inverses \cite{Reu02}.
We also show that the accuracy of
sparse inverse approximations increases when more matrix
elements are included.

\paragraph{Idea}
To understand how sparse inverse approximations work,
we first consider a representation of the determinant based on
minors of the inverse \cite[\S 0.8.4]{HoJ85}.
If $M$ is Hermitian positive-definite of order $n$, 
and $M_i$  is the leading principal submatrix of order $i$ of $M$,
then \cite[\S 0.8.4]{HoJ85} $\det(M)=\det(M_{n-1})/\sigma_n$,
where $\sigma_n\equiv (M^{-1})_{nn}$ is the trailing diagonal element
of $M^{-1}$. Using this expression recursively for $\det(M_{n-1})$
gives
$$\det(M)=\prod_{i=1}^{n}{\frac{1}{\sigma_i}},\qquad
\mathrm{where}\quad \sigma_i=(M_i^{-1})_{ii}.$$

Determinant approximations based on sparse approximate inverses
replace \textit{leading} principal submatrices $M_i$ 
by just principal submatrices $S_i$. 
Specifically \cite[\S 3.2]{Reu02}, let $M$ be Hermitian positive-definite,
and let $S_i$ be a principal submatrix of $M_i$, 
such that $S_i$ includes at least row $i$ and column $i$ of $M$.
The two extreme cases are $S_i=m_{ii}$ and $S_i=M_i$.
In any case, $m_{ii}$ is the trailing diagonal element of $S_i$,
i.e. $S_{n_i,n_i}=m_{ii}$, where $n_i$ is the order of $S_i$,
$1\leq n_i\leq i$.
Let $\sigma_i$ be the trailing diagonal element of $S_i^{-1}$, 
i.e. $\sigma_i=(S_i^{-1})_{n_i,n_i}$. In particular $\sigma_1=m_{11}^{-1}$.
Given $n$ such submatrices $S_i$, $1\leq i\leq n$,
the sparse inverse approximation of $\det(M)$ is defined as
\cite[Algorithm 3.3]{Reu02}.
\begin{eqnarray}\label{e_si}
\sigma\equiv \prod_{i=1}^n{\frac{1}{\sigma_i}}.
\end{eqnarray}
The sparse approximate inverse method performs
Cholesky decompositions $S_i=L_iL_i^*$, where $L_i$ is lower triangular,
$2\leq i\leq n$, and computes $1/\sigma_i=((L_{i})_{n_i,n_i})^2$.

\paragraph{Monotonicity}
We show monotonicity of the sparse inverse approximations 
in the following sense: If the dimensions of the submatrices $S_i$ 
are increased then the determinant  approximations can only become better.

\begin{lemma}\label{l_invdiag}
If 
$$M=\bordermatrix{&m&k\cr m&A&B\cr k&B^* &S}$$
is Hermitian positive-definite then
$(S^{-1})_{ii}\leq (M^{-1})_{m+i,m+i}$, $1\leq i\leq k$.
\end{lemma}

\begin{proof}
The proof follows from \cite[(4)]{Cot74} and the
Shermann-Morrison formula \cite[(2.1.4)]{GovL96}.
\end{proof}

Lemma \ref{l_invdiag} implies the following lower and upper bounds for
sparse inverse approximations; the lower bound was already derived in
\cite[(3.25)]{Reu02}.

\begin{corollary}\label{c_si}
If $M$ is Hermitian positive-definite and $\sigma$ is a sparse
inverse approximation in (\ref{e_si}) then
$$\det(M)\leq \sigma\leq \prod_{i}{m_{ii}}.$$
\end{corollary}

Corollary \ref{c_si} implies that the product of diagonal elements
cannot approximate the determinant more accurately than a sparse
inverse approximation. Another consequence of 
Lemma \ref{l_invdiag} is the monotonicity of the sparse inverse
approximation in the following sense: If a principal submatrix ${\hat S}_j$ 
is replaced by a larger principal submatrix  $S_j$ then the 
determinant approximation can only become better.

\begin{theorem}\label{t_si}
Let $M$ be Hermitian positive-definite of order $n$.
If for some $1<j\leq n$, $S_j$ is a principal submatrix of $M_j$, 
and in turn ${\hat S}_j$ is a principal submatrix of $S_j$ then
$$\det(M)\leq \prod_{i=1}^n{1\over\sigma_i}\leq
{1\over {\hat\sigma}_j}\prod_{i=1,i\neq j}^n{1\over\sigma_i}$$
where ${\hat\sigma}_j$ is the trailing diagonal
element of ${\hat S}_j^{-1}$.
\end{theorem}

The next example of block diagonal matrices illustrates that
sparse inverse approximations can be inaccurate, even when sparsity
is exploited to full extent.

\paragraph{Block-Diagonal Matrices}
Let 
$$M=\begin{pmatrix}T_3 && \\ &\ddots& \\&&T_3\end{pmatrix}$$
be a block diagonal matrix of order $n=3k$ with $n/3$ diagonal blocks 
$$T_3\equiv\begin{pmatrix}3/2&-1 &\cr-1&3/2& -1\cr &-1&3/2\end{pmatrix}.$$
The obvious block diagonal approximation (\ref{e_part}) with $k=n/3$ gives the
exact determinant $\det(\md)=\det(M)=\det(T_3)^{n/3}=(3/8)^{n/3}$.
For the sparse inverse approximation (\ref{e_si}) we choose the submatrices
$$S_{(i-1)(n/3)+1}=3/2,\qquad 
S_{(i-1)(n/3)+2}=\begin{pmatrix}3/2 &-1\cr -1& 3/2\end{pmatrix}=S_{i(n/3)},
\qquad 1\leq i\leq n/3.$$
The sparse inverse approximation of $\det(T_3)$ is $\det(T_3)+2/3$.
It has no accurate digit because the relative error is $16/9$.
The sparse inverse approximation of $\det(M)$ is
$\sigma=\left(\det(T_3)+2/3)\right)^{n/3}$.
For instance, when $n=300$ then $\det(M)\approx 4\cdot10^{17}$
while the sparse inverse approximation gives $\sigma\approx 4\cdot10^{33}$.

\paragraph{Tridiagonal Toeplitz Matrices}
A block diagonal approximation can be more accurate than a
sparse inverse approximation if the dimension of the blocks is
larger than 1. Let
$$T_n=\begin{pmatrix}2 & -1  &\cr
             -1  & 2 & \ddots \cr
                    &\ddots&\ddots&-1\cr
                    &      & -1&2\end{pmatrix}$$
be of order $n$; then $\det(T_n)=n+1$.
In the sparse inverse approximation (\ref{e_si}) we fully exploit sparsity by choosing
$S_1=2$ and $S_i=T_2$, $2\leq i\leq n$; hence the approximation
is $\sigma=2\>(3/2)^{n-1}$.
When $\md$ in (\ref{e_part}) consists of $k$ equally sized blocks
of dimension $n/k$ then
$\det(\md)=\left(\det(T_{n/k})\right)^k=\left((n/k)+1\right)^k$.
For a block size $n/k\geq 4$, $\det(\md)\leq \sigma$, so
the block diagonal approximation is more accurate than
the sparse inverse approximation.

\paragraph{2-D Laplacian}
We show that for this matrix
both, the block-diagonal and the sparse inverse
approximations are accurate to at most one digit.
The coefficient matrix from the centered finite difference discretization of 
Poisson's equation is a Hermitian positive-definite block
tridiagonal matrix \cite[9.1.1]{Gre97}
$$M=\begin{pmatrix}T_m &-I_m\cr -I_m&T_m&\ddots\cr 
                 &\ddots&\ddots &-I_m\cr &&-I_m&T_m\end{pmatrix},\qquad
\mathrm{where}\quad
T_m=\begin{pmatrix}4 &-1\cr -1&4&\ddots\cr 
                 &\ddots&\ddots &-1\cr &&-1&4\end{pmatrix}.$$
Here $T_m$ is of order $m$, and $M$ is of order $n=m^2$ 
(note that the matrix considered in \cite[\S 5]{Reu02} equals $(n+1)^2M$).
The exact determinant is \cite[Theorem 9.1.2]{Gre97}
$$\det(M)=\prod_{i,j=1}^m{4\left(\sin{\left({i\pi\over 2(m+1)}\right)}^2+
\sin{\left({j\pi\over 2(m+1)}\right)}^2\right)}.$$
We compute the logarithm of this expression and compare it to the
approximations.
A block diagonal approximation (\ref{e_part}) with $k=m$ gives
$\det(\md)=\det(T_m)^m$, where \cite[\S 28.5]{Hig02}
$$\det(T_m)=\prod_{i=1}^m{\left(4+2\cos{i\pi\over m+1}\right)}.$$
If the matrices in the sparse inverse approximation (\ref{e_si}) are
$$S_1=4,\qquad S_i=\begin{pmatrix}4&-1\cr -1&4\end{pmatrix},\quad 
2\leq i\leq m+1,\qquad 
S_j=\begin{pmatrix}4&0&-1\cr 0&4&-1\cr -1&-1&4\end{pmatrix},\quad 
m+2\leq j\leq n,$$
then $1/\sigma_1=4$, $1/\sigma_i=15/4$, $2\leq i\leq m+1$, and
$1/\sigma_j=7/2$, $m+2\leq j\leq n$. Thus the sparse
inverse approximation is $\sigma=4(15/4)^m(7/2)^{n-m-1}$.

Table \ref{table_2} lists errors for 
the block diagonal and sparse inverse approximations for 
$n=900$, $n=10000$ and $n=40000$. Columns 3 and 4 represent the 
relative errors 
$$|\ln(\det(\md))-\ln(\det(M))|/|\ln(\det(M))|\qquad\mathrm{and}\qquad
|\ln(\sigma)-\ln(\det(M))|/|\ln(\det(M))|,$$ 
while columns 5 and 6 represent the relative errors
$$|\det(\md)^{1/n}-\det(M)^{1/n}|/|\det(M)^{1/n}|\qquad\mathrm{and}\qquad
|\sigma^{1/n}-\det(M)^{1/n}|/|\det(M)^{1/n}|.$$ 
We include the last two errors to allow a comparison with the 
approximation of $\det\left((n+1)^2M\right)^{1/n}$ in \cite[Table 5.1]{Reu02}.
The table shows that all relative errors lie between $0.06$ and $0.2$.
Hence both approximations, block diagonal and sparse inverse,
are accurate to at most one significant digit.
To estimate the tightness of the  bound 
$$|\ln(\det(\md))-\ln{\det(M)}|\leq (-n\>|\ln(1-\rho)|)\>\rho$$
in Theorem \ref{t_3},
consider the case $n=900$. Here $\rho(\md^{-1}\moff)\approx .9898$
and $|\ln(1-\rho)|\approx 4.5845$.
The true error is 
$$|\ln(\det(\md))-\ln{\det(M)}|\approx 122.4966\approx 26\ln(1-\rho)\>\rho.$$
The matrix $\md^{-1}\moff$ has 26 eigenvalues with magnitude at least $0.9$.
Thus the pessimism of the bound comes from the factor $n$.

\begin{table}
\begin{center}
\begin{tabular}{|c|l|l|l|l|l|}\hline
$n$&$\ln(\det(M))$&rel. error& rel. error& rel error&rel. error\\
&&in $\ln(\md)$& in $\ln(\sigma)$&in $\md^{1/n}$& in $\sigma^{1/n}$\\
\hline
900&    1.0650e+03&0.1150&    0.0607&    0.1458&    0.0745\\
10000&    1.1717e+04&   0.1246&    0.0698&    0.1572&    0.0852\\
40000&   4.6761e+04&    0.1269&    0.0719&    0.1599&    0.0877\\
\hline
\end{tabular}
$$\qquad$$
\end{center}
\caption{Errors in the block diagonal approximation $\md$ and 
the sparse inverse approximation $\sigma$ 
for the Laplacian.}\label{table_2}
\end{table}

\section{Application to Neutron Matter Simulations}\label{s_app}
In \cite{LeeI03} we consider the quantum simulation of
nuclear matter on a lattice, and in particular how to calculate the
contribution of nucleon-nucleon-hole loops at non-zero nucleon
density. The resulting method, called zone determinant expansion,
is based on the sequence of approximations in Theorem \ref{t_3}.
Here we illustrate that 3 iterations of the zone determinant
expansion give an approximation accurate to 3 digits, and that the
method uses less space than a determinant computation based on Gaussian
elimination (with partial or complete pivoting).

In \cite{LeeI03} we derive a particle interaction matrix $M$
whose determinant $\det(M)$ is not positive, and complex in general.
Hence stochastic methods such as
hybrid Monte Carlo methods \cite{DuKPR87,GoLTRS87,ScSS86}
do not give the correct sign or phase of $\det(M)$.
This was the motivation for approximating
$\ln(\det(M))$ via a zone determinant expansion, i.e. Theorem \ref{t_3}.
Below we discuss the structure of $M$ and a physically appropriate
zone determinant expansion.

The particle interactions are considered on a 4-dimensional lattice
(3 dimensions for space and one for time).
Let the dimensions of the lattice be $L\times L\times L\times L_t$,
where $L_t$ represents the time direction.
Also let the number of particles per lattice point be $s$.
Then the interaction matrix $M$ has dimension $n\times n$
where $n=L^3L_ts$. We partition the lattice into separate spatial 
zones (or cubes) of dimension $m\times m \times m$ (constraints on $m$
are discussed in \cite{LeeI03}).
Therefore particle interactions between any two zones are
represented by matrix blocks of dimension $m^3L_ts$.
As a consequence, it makes sense to approximate $\det(M)$ by the
product of principal minors associated 
with particle interactions inside spatial zones.
Without loss of generality we assume that the lattice points
are ordered such that the submatrix $M_{ij}$ of order $m^3L_ts$
represents particle interactions between zones $i$ and $j$.
With $k\equiv (L/m)^3$ this gives the partitioning
$M=\md+\moff$ in (\ref{e_part}), where $\md$ represents
particle interactions in the zone interiors,
while $\moff$ represents interactions among different zones.
In \cite{LeeI03} we explain that the spectral radius 
$\rho\equiv\rho(\md^{-1}\moff)$ can be reduced
by increasing the dimension $m$ of the spatial zones.

\begin{figure}
\begin{center}
\resizebox{3in}{!}
{\includegraphics*{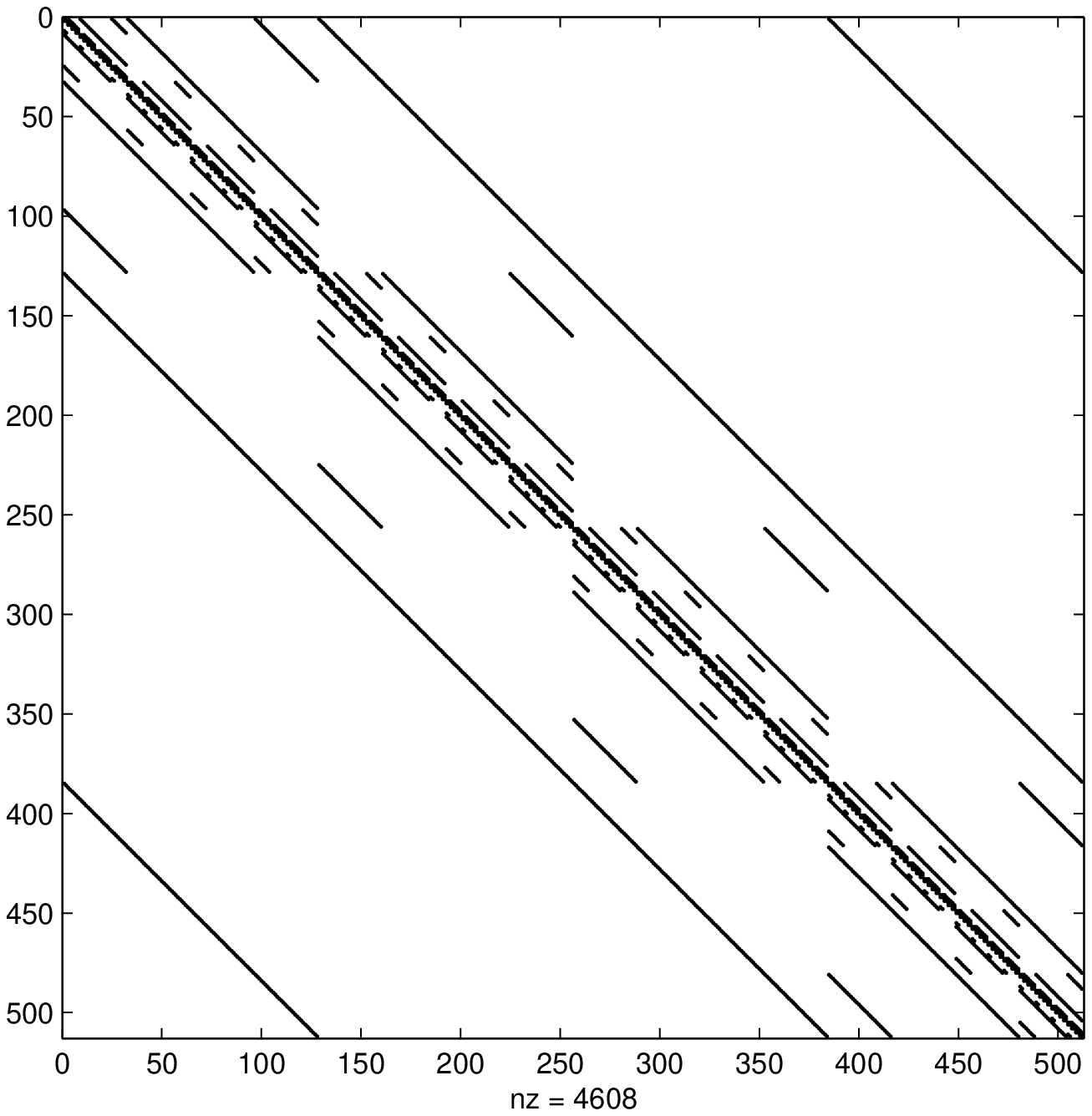}}
\end{center}
\caption{Sparsity structure of the interaction matrix $M$.}\label{f_spy}
\end{figure}

We illustrate the zone expansion on a small lattice simulation,
where we can compare the approximations to the exact determinant.
Specifically we consider the interactions between  neutrons and neutral 
pions, on a $4^3\times 4$ grid. The order of the interaction matrix $M$ is
$4^3\times 4\times 2=512$. Its properties are listed in Table \ref{table_3}.

\begin{table}
\begin{tabular}{|l|l|l|}
\hline
order & $n=512$ & \\
number of non-zeros & $9n$ & see Figure \ref{f_spy}\\
structure & complex non-Hermitian & see Figure  \ref{f_blockspy}\\
norm & $\|M\|_F\approx 49.5$ & \\
condition number & $\|M\|_1\|M^{-1}\|_1\approx 177$ & 
   condest$(\cdot)$  in MATLAB 6\\
non-normality & $\|M^*M-MM^*\|_F\approx 57$ &\\
eigenvalues & complex & see Figure \ref{f_evm} \\
determinant   & $\det(M)=8.5361\cdot 10^{65}+1.4168\cdot 10^{64}\imath$ 
                     & $\det(\cdot)$ in Matlab 6\\
&$\ln(\det(M))=151.81+0.016599\imath$ &\\
\hline
\end{tabular}
\caption{Properties of the interaction matrix $M$.}\label{table_3}
\end{table}

\noindent
\begin{figure}
\begin{center}
\resizebox{2in}{!}
{\includegraphics*{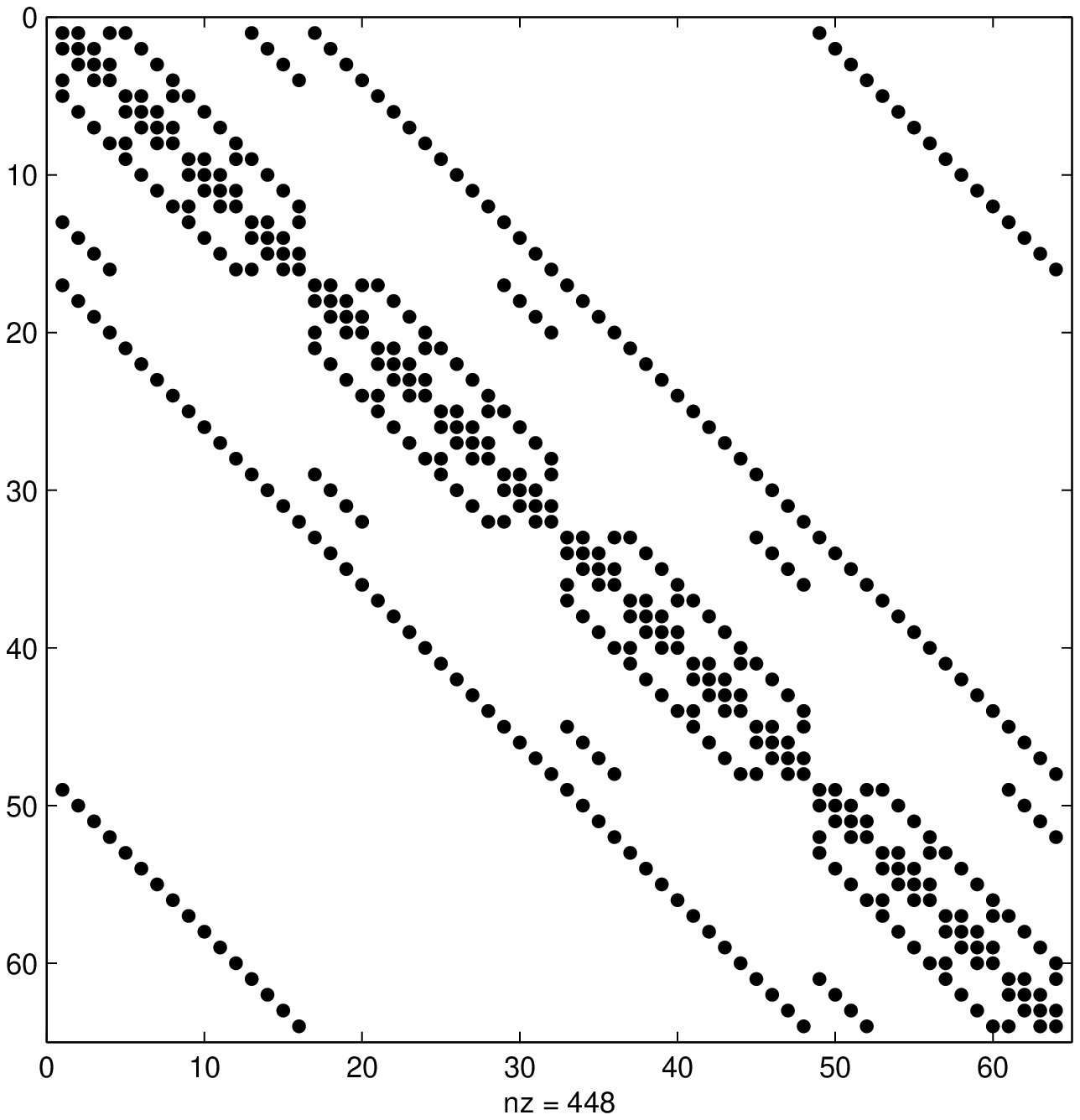}}
$\qquad\qquad\qquad$
\resizebox{1in}{!}
{\includegraphics*{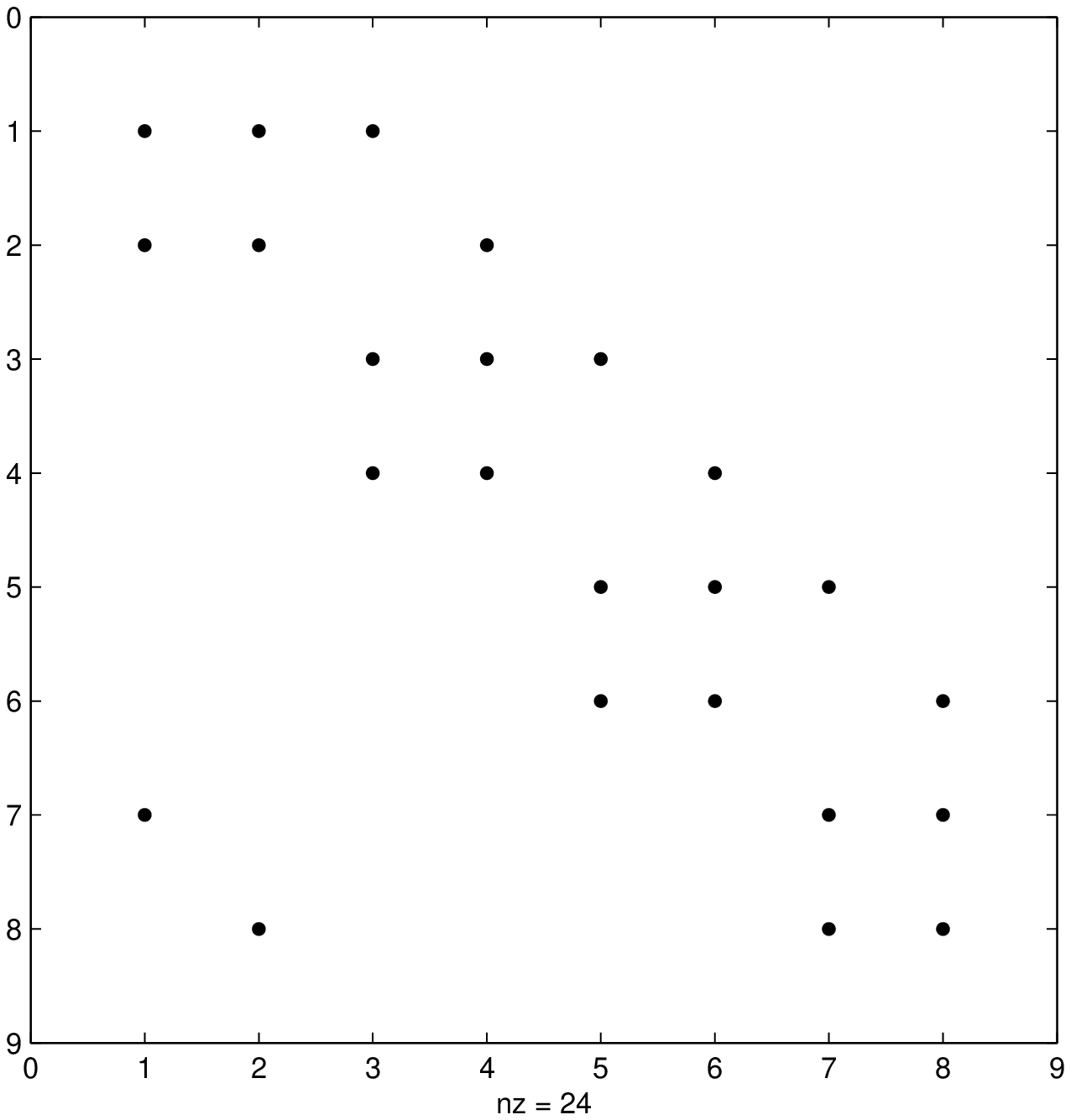}}
\end{center}
\caption{Non-zero $8\times 8$ blocks in the interaction matrix $M$, and 
sparsity structure of a single $8\times 8$ diagonal block.}\label{f_blockspy}
\end{figure}

In the context of the particular application in \cite{LeeI03},
we can partition the lattice into zones with dimension $m=1$.
The resulting partitioning has blocks 
$M_{ij}\equiv M_{8(i-1)+1:8i,8(j-1)+1:j}$, $1\leq i,j\leq 64$,
of dimension $4\times 2=8$. Thus $k=64$ in the block diagonal
approximation (\ref{e_part}).
Figure \ref{f_blockspy} shows the distribution of the 448 blocks 
with non-zero elements.
Each diagonal block $M_{ii}$ contains 24 non-zero elements, 
its sparsity structure is  shown in Figure \ref{f_blockspy}.

\begin{figure}
\begin{center}
\resizebox{3in}{!}
{\includegraphics*{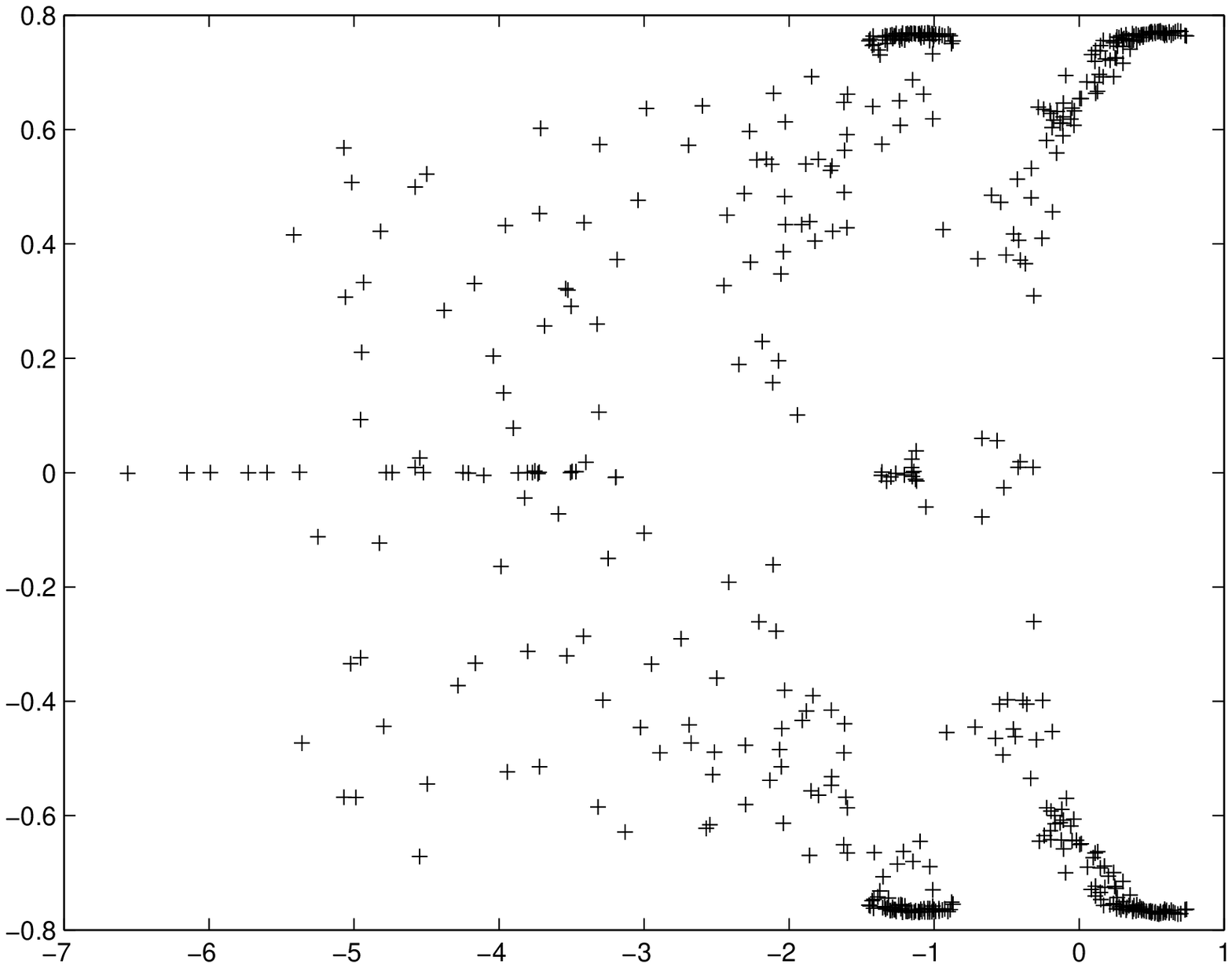}}
\end{center}
\caption{Eigenvalue distribution of the interaction matrix $M$.}\label{f_evm}
\end{figure}

\begin{figure}
\begin{center}
\resizebox{1.5in}{!}
{\includegraphics*{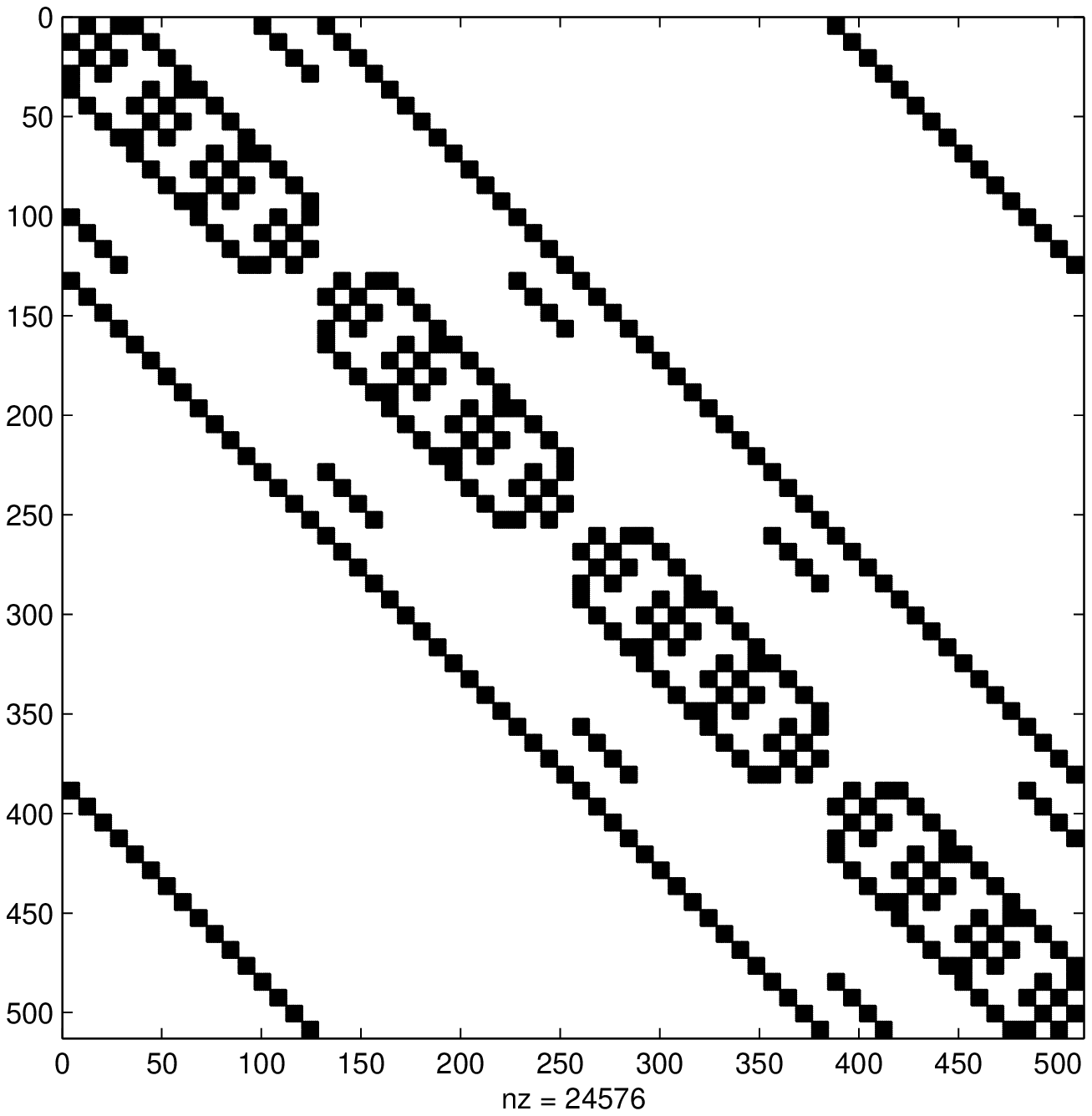}}
$\qquad$
\resizebox{1.5in}{!}
{\includegraphics*{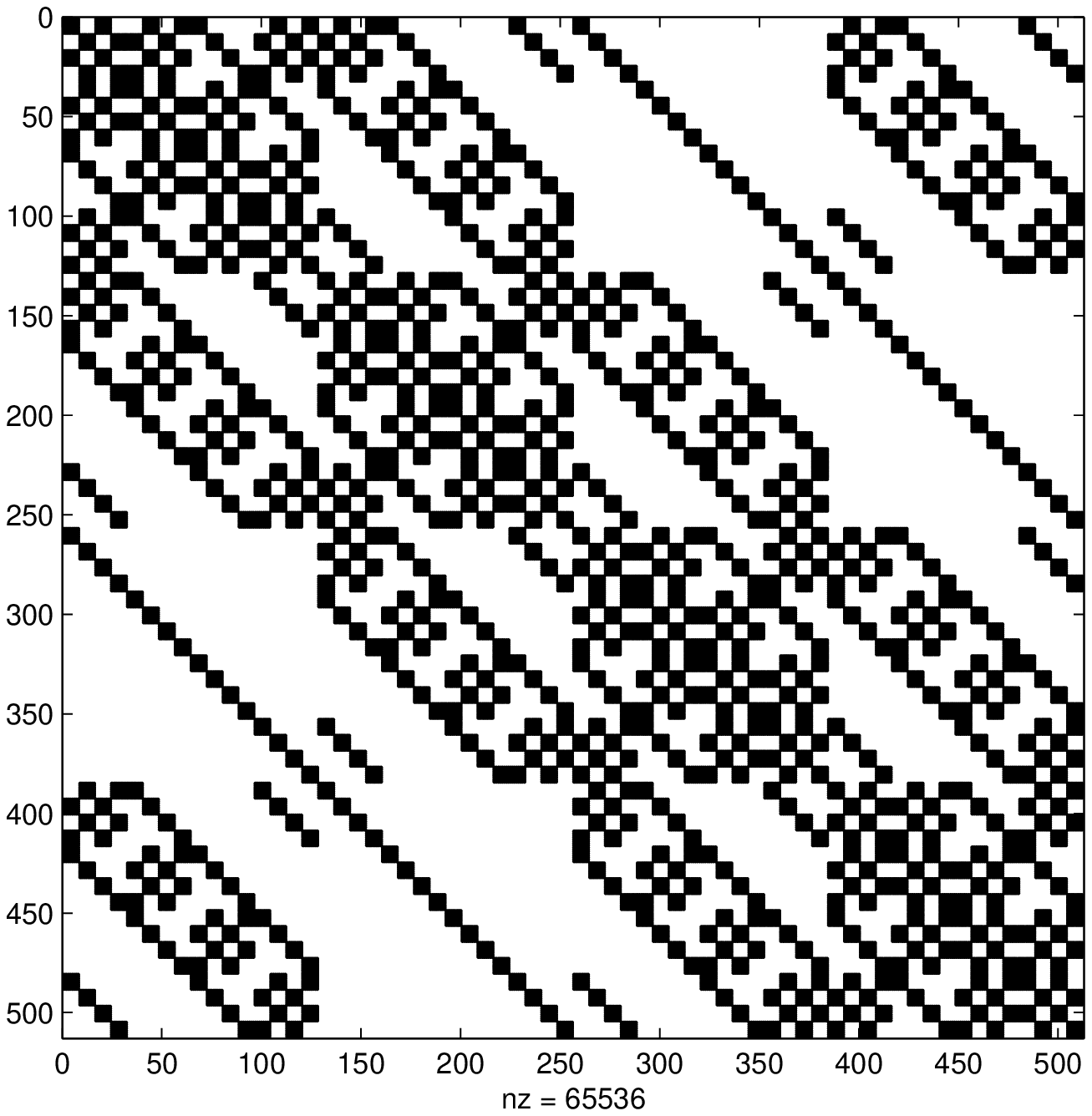}}
\end{center}
\caption{Sparsity structure of the matrices $\md^{-1}\moff$ and
$(\md^{-1}\moff)^2$.}\label{f_zone}
\end{figure}

\begin{figure}
\begin{center}
\resizebox{1in}{!}
{\includegraphics*{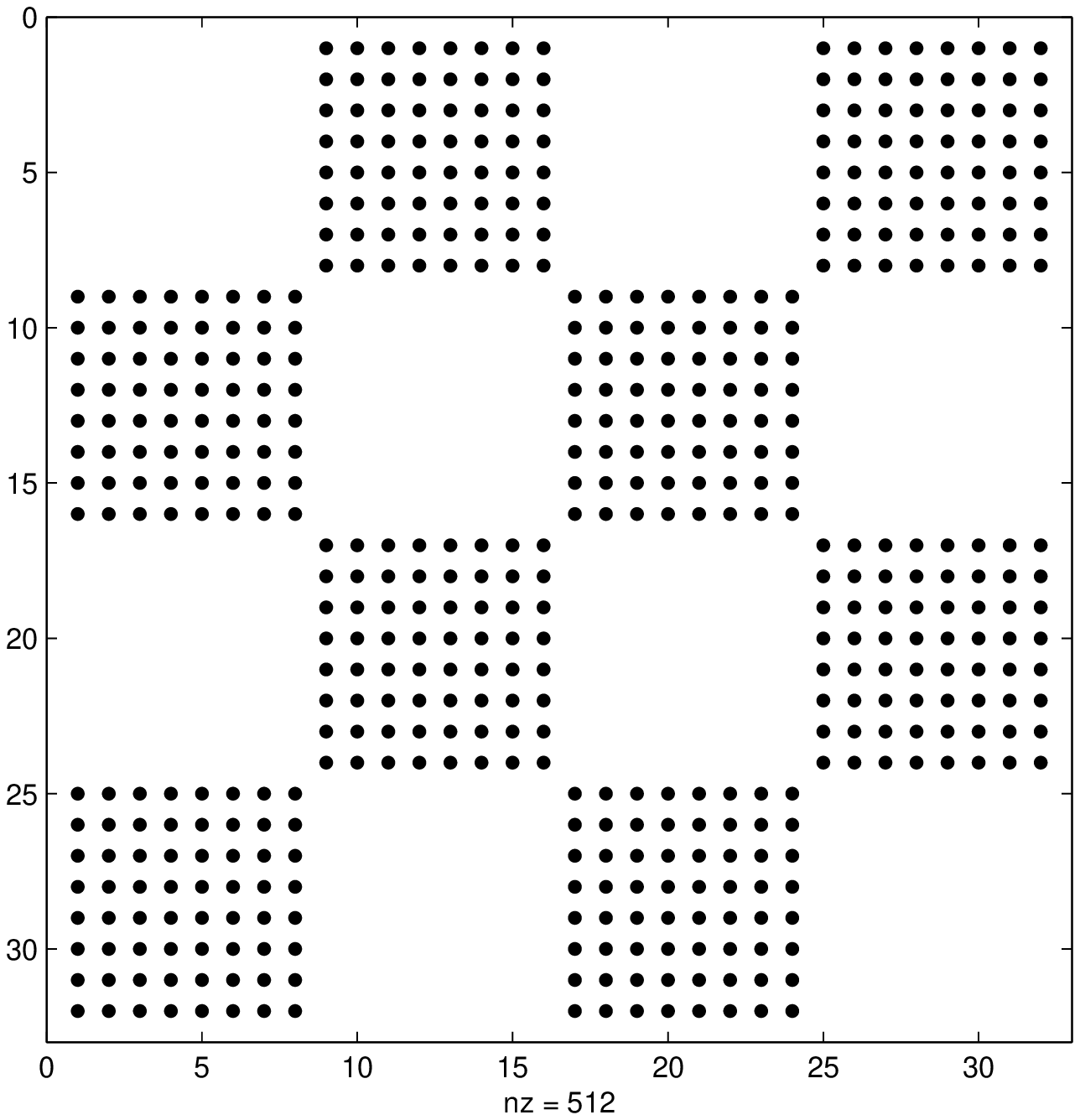}}
\end{center}
\caption{Sparsity structure of the leading principal submatrix of order 32
of  $\md^{-1}\moff$.}\label{f_checker}
\end{figure}

The zone partitioning is bipartite, i.e.
$M_{ij}=0$ for $i$ and $j$ both even or both odd, and $i\neq j$, 
$1\leq i,j\leq k$. 
Therefore $\moff$ is an odd checkerboard matrix. 
Figure \ref{f_checker} illustrates this checkerboard pattern in the 
leading principal submatrix of order 32 of $\md^{-1}\moff$.
The sparsity structures of the matrices
$\md^{-1}\moff$ and  $(\md^{-1}\moff)^2$
is shown in Figure \ref{f_zone}. Because of the checkerboard structure
Theorem \ref{t_checker} implies 
$\trace((\md^{-1}\moff)^p)=0$ for odd $p$, and 
$\delta_{p-1}=\delta_{p}$. Table \ref{table_1} therefore contains
only approximations of even order.

\begin{table}
\begin{center}
\begin{tabular}{|c|l|l|l|l|l|l|}\hline
$j$ &abs. error &abs. error &abs. error &$\rho^j$& rel. error &rel. error \\
&in $\Re(\delta_j)$ &in $\Im(\delta_j)$ &in $\delta_j$ &&in $\delta_j$&
in $\Delta_j$ \\
\hline
0&5.1000 & 0.0017&5.1000&& 0.0348&163.0282\\
2& 0.4817& 0.0025&0.4817&0.4374&0.0032&0.3823\\
4&  0.0909&0.0016&0.0909&0.1913&0.0006&0.0951\\
6& 0.0225&0.0008&0.0226&0.0837&0.0001&0.0223\\
8&    0.00665&  0.0003&0.0066& 0.0366&0.00004& 0.0066\\
\hline
\end{tabular}
$$\qquad$$
\end{center}
\caption{Errors in the approximations $\delta_j$ and $\Delta_j$
for the interaction matrix $M$.}\label{table_1}
\end{table}

Table \ref{table_1} shows errors in the approximations $\delta_j$
and $\Delta_j$ for approximations up to order 8.
Columns 2, 3 and 4 represent the absolute errors
$$|\Re(\ln(\det(M)))-\Re(\delta_j)|, \qquad 
|\Im(\ln(\det(M)))-\Im(\delta_j)|,\qquad
|\ln(\det(M))-\delta_j|.$$ 
Columns 6 and 7 represent the relative errors
$$|\ln(\det(M))-\delta_j|/|\delta_j|\qquad 
\mathrm{and}\qquad  |\det(M)-\Delta_j|/|\Delta_j|.$$
The spectral radius $\rho\equiv\rho(\md^{-1}\moff)\approx .6613$, 
and the constant
in the error bounds of Theorems \ref{t_2} and \ref{t_3} is $c\approx 554$.

Table \ref{table_1} illustrates that
$|\ln(\det(M))-\delta_j|\approx\rho^j$,
i.e. the absolute errors in  the logarithm are almost proportional
to the powers of the spectral radius of $\md^{-1}\moff$.
In this case the constant $c$ is too pessimistic, because
many eigenvalues of $\md^{-1}\moff$ have magnitude much less than $\rho$.
For instance, 160 eigenvalues of $\md^{-1}\moff$ have magnitude $10^{-15}$.
The imaginary parts of the logarithms appear to converge faster
than the real parts.
The block diagonal approximation $\delta_0\equiv \ln(\det(\md))$ 
for $\ln(\det(M))$ has
an accuracy of 2 digits. Two more iterations give an approximation 
$\delta_2$  that is accurate to 3 digits.

\begin{figure}
\begin{center}
\resizebox{1.5in}{!}
{\includegraphics*{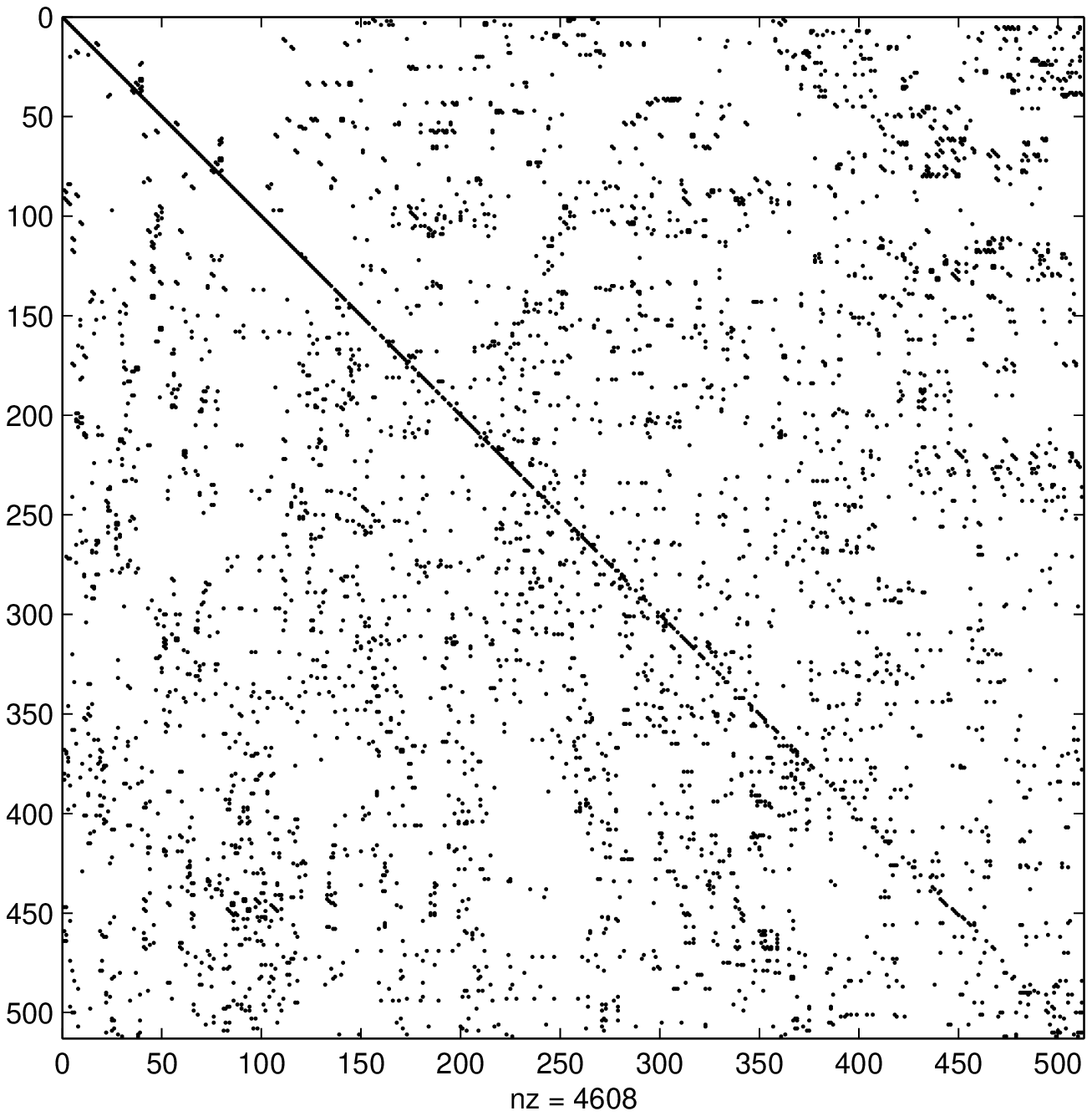}}
$\qquad$
\resizebox{1.5in}{!}
{\includegraphics*{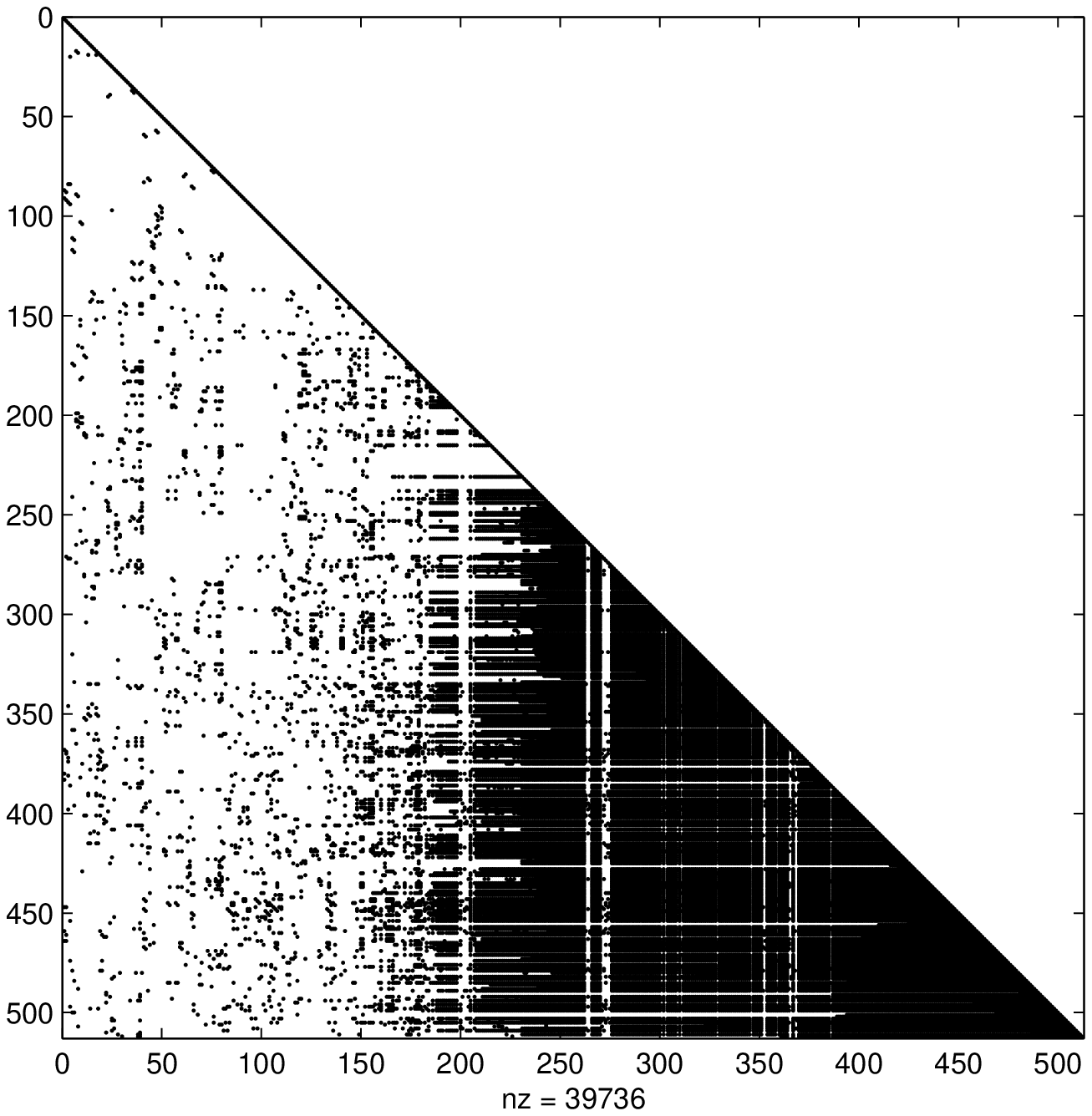}}
$\quad$
\resizebox{1.5in}{!}
{\includegraphics*{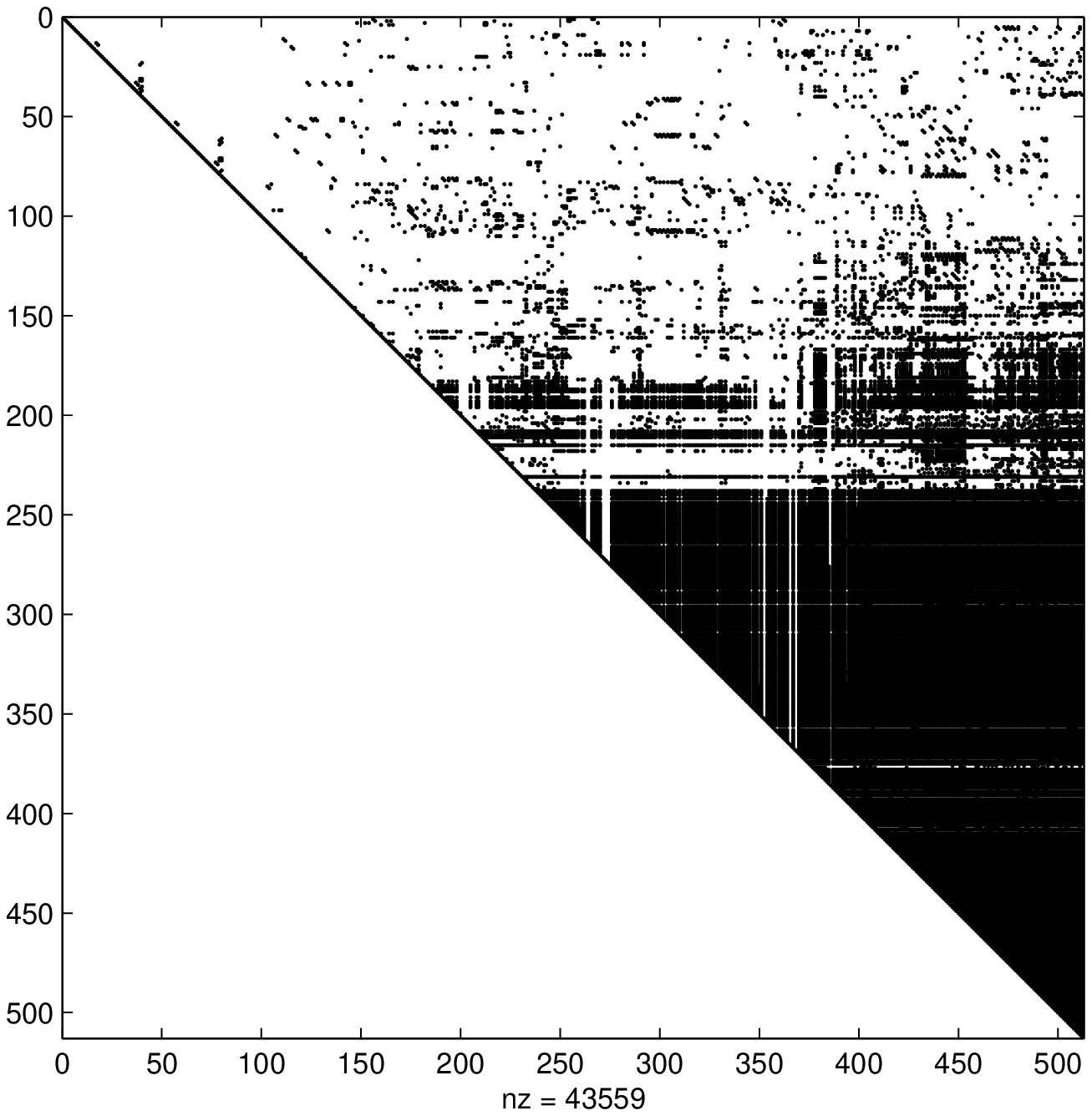}}
\end{center}
\caption{Sparsity structure of the matrices $PMQ$, $L$ and $U$ from the
LU decomposition (with complete pivoting) of $M$.}\label{f_lu}
\end{figure}

We briefly compare the computation of $\delta_0$ and $\delta_2$ 
to a determinant computation by Gaussian elimination of $M$.
Gaussian elimination with complete pivoting gives $PMQ=LU$, where 
$P$ and $Q$ are permutation matrices, $L$ is unit lower triangular and $U$ is 
upper triangular.
Figure \ref{f_lu}, which shows the sparsity
structure of the matrices $PMQ$, $L$ and $U$, illustrates that
Gaussian elimination with complete pivoting completely destroys the
sparsity structure of $M$. The matrices $L$ and $U$ together have
about $162n$ non-zeros, compared to $9n$ in $M$.
In contrast, the determinant expansion
requires no significant additional space for $\delta_0$; and $48n$
non-zeros for $\md^{-1}\moff$ and $n$ non-zeros for 
the trace of $(\md^{-1}\moff)^2$. That's $(48+1)n=49n$ non-zeros,
about one third of the non-zeros produced by Gaussian elimination with
complete pivoting. Gaussian elimination with partial pivoting
essentially preserves the sparsity structure of $M$ but produces
$342n$ non-zeros.

\subsection*{Acknowledgements}
We thank Gene Golub, Nick Higham, Volker Mehrmann, and Gerard Meurant for
helpful discussions.

\end{document}